\documentclass[11pt,reqno]{amsart}
\usepackage[T1]{fontenc}
\usepackage{lmodern}
\usepackage{amsmath,amssymb,mathtools}
\usepackage{microtype}
\usepackage{needspace}
\usepackage[hidelinks]{hyperref}

\numberwithin{equation}{section}
\newtheorem{theorem}{Theorem}[section]
\newtheorem{lemma}[theorem]{Lemma}
\newtheorem{corollary}[theorem]{Corollary}

\newcommand{\R}{\mathbb R}
\newcommand{\C}{\mathbb C}
\newcommand{\eqnum}[1]{\refstepcounter{equation}\label{#1}\eqno{\textup{(\theequation)}}}

\title[Polynomial approximation with monotone weights]{A counterexample to polynomial approximation for monotone weights}

\author[N. Alpay]{Natanael Alpay}
\address{(NA)
Department of Mathematics\\
University of California, Irvine,
Irvine, CA 92697\\
USA}
\email{nalpay@uci.edu}

\author[P. Ivanisvili]{Paata Ivanisvili}
\address{(PI)
Department of Mathematics\\
University of California, Irvine,
Irvine, CA 92697\\
USA}
\email{pivanisv@uci.edu}

\date{}

\begin{document}

\begin{abstract}
We give a counterexample to a question of Borichev on polynomial
approximation with monotone weights. We construct an even
continuous weight $W$, strictly increasing on $[0,\infty)$,
with superpolynomial growth and finite logarithmic integral,
and an entire function $F$ of zero exponential type such that
$F/W$ vanishes at infinity and belongs to $L^p(\R)$ for every
$1\leq p<\infty$. Nevertheless, $F$ has positive distance from
the polynomials in the weighted uniform norm and in every
corresponding weighted $L^p$ norm.
\end{abstract}

\maketitle

\section{Introduction}\label{sec:intro}

We consider the closure of polynomials in weighted spaces of
functions on the real line. Let $W:\R\to[1,\infty)$ be an even
continuous weight satisfying
\begin{equation}\label{eq:growth}
\lim_{|x|\to\infty}\frac{\log W(x)}{\log|x|}=\infty.
\end{equation}
We work with complex-valued functions and write
\[
C_W^0
=
\left\{
f\in C(\R):
\frac{f(x)}{W(x)}\longrightarrow0
\text{ as }|x|\to\infty
\right\},
\qquad
\|f\|_W=\sup_{x\in\R}\frac{|f(x)|}{W(x)}.
\]
For $1\leq p<\infty$, set
\[
L_W^p=L^p(\R,W^{-p}\,dx),
\qquad
\|f\|_{p,W}
=
\left(
\int_\R\frac{|f(x)|^p}{W(x)^p}\,dx
\right)^{1/p}.
\]
Condition \eqref{eq:growth} ensures that every polynomial belongs
to each of these spaces.

An entire function $F$ has \emph{zero exponential type} if, for
every $\varepsilon>0$, there is a constant $C_\varepsilon>0$ such
that
\[
|F(z)|\leq C_\varepsilon e^{\varepsilon|z|},
\qquad z\in\C.
\]
We are concerned with weights for which
\begin{equation}\label{eq:logint}
\int_\R\frac{\log W(x)}{1+x^2}\,dx<\infty.
\end{equation}
The approximation problem is to determine whether every entire
function of zero exponential type belonging to a given weighted
space lies in its polynomial closure.

A weight is called log-convex here if
$t\mapsto\log W(e^t)$ is convex on $\R$.
For log-convex weights satisfying \eqref{eq:growth} and
\eqref{eq:logint}, Koosis \cite[Theorems~I and~IV]{Koosis1964} obtained
the approximation result in $C_W^0$ under the additional
condition that, for every $\theta>1$, there is a constant
$C_\theta>0$ such that
\[
x^2W(x)\leq C_\theta W(\theta x),
\qquad x\geq0.
\]
Borichev \cite[Theorems~2.1 and~3.1]{Borichev2001} removed this
additional condition
and proved that the polynomial closure in $C_W^0$ and in every
$L_W^p$, $1\leq p<\infty$, consists exactly of the restrictions
of entire functions of zero exponential type belonging to the
respective space. In \cite[Introduction, p.~831]{Borichev2001},
he asked whether log-convexity could be replaced by monotonicity
of $W$ on $[0,\infty)$, while retaining \eqref{eq:growth} and
\eqref{eq:logint}.

Counterexamples without monotonicity were already known; see,
for instance, Koosis \cite[Theorem~VI]{Koosis1964}.
Koosis also remarked, without giving the details, that his second
construction can be made monotone
\cite[remark following Theorem~VII, p.~407]{Koosis1964}.
For that example, finite sums of exponentials with frequencies in
$[-a,a]$ are dense in the whole weighted continuous space for every
$a>0$, whereas the polynomials are not. The logarithmic integral
must therefore be infinite, by
\cite[lemma preceding Theorem~II, p.~391]{Koosis1964}, so this
example lies outside the setting of \eqref{eq:logint}.

A different phenomenon occurs in an earlier construction of
Borichev \cite[Theorem~1 and the concluding remark]{Borichev2000}.
He constructed an even continuous weight $V$, nondecreasing on
$[0,\infty)$, for which polynomials are dense in $C_V^0$, although
\[
\int_\R\frac{\log V(x)}{1+x^2}\,dx=\infty,
\qquad
\int_\R\frac{\log V_{\mathrm{lc}}(x)}{1+x^2}\,dx<\infty,
\]
where $V_{\mathrm{lc}}\leq V$ is the largest log-convex minorant
of $V$. Thus the finite logarithmic integral concerns the minorant,
not the weight governing approximation. The question considered
here requires \eqref{eq:logint} for $W$ itself.

We give a negative answer to Borichev's question. We construct
an even continuous weight $W$, strictly increasing on
$[0,\infty)$ and satisfying \eqref{eq:growth} and
\eqref{eq:logint}, together with an entire function $F$ of zero
exponential type that belongs to $C_W^0$ and every $L_W^p$,
$1\leq p<\infty$, but has positive distance from the polynomials
in each of these spaces. The same weight and function therefore
give a counterexample simultaneously for the weighted uniform
norm and all the finite weighted $L^p$ norms.

\Needspace{12\baselineskip}
\section{Main results}\label{sec:results}

\begin{theorem}\label{thm:main}
There exist an even continuous weight $W:\R\to[1,\infty)$, strictly
increasing on $[0,\infty)$ and satisfying \eqref{eq:growth} and
\eqref{eq:logint}, and an entire function $F$ of zero exponential type
such that
$$
F\in C_W^0\cap\bigcap_{1\leq p<\infty}L_W^p,
$$
but
$$
\inf_{P\in\C[x]}\|F-P\|_W>0,
\qquad
\inf_{P\in\C[x]}\|F-P\|_{p,W}>0
\quad(1\leq p<\infty).
\eqnum{eq:mainseparation}
$$
Moreover, for every $a>0$, $F$ belongs, in each of these spaces, to the
closure of
$$
\operatorname{span}\{e^{itx}:t\in\R,\ |t|\leq a\}.
\eqnum{eq:smallfreqspace}
$$
\end{theorem}

The last part of Theorem~\ref{thm:main} follows from the weighted
approximation characterization stated in
\cite[Introduction, p.~831]{Borichev2001}, where it is attributed
to de~Branges. Under \eqref{eq:growth} and \eqref{eq:logint}, the
entire functions of zero exponential type in each weighted space
form the intersection of the closures of the exponential spans
in \eqref{eq:smallfreqspace}.

\begin{corollary}\label{cor:strictinclusion}
For the weight in Theorem~\ref{thm:main}, let $X$ be either $C_W^0$ or
$L_W^p$, $1\leq p<\infty$. Then
$$
\overline{\C[x]}^{\,X}
\subsetneq
\bigcap_{a>0}
\overline{\operatorname{span}\{e^{itx}:t\in\R,\ |t|\leq a\}}^{\,X}.
$$
\end{corollary}

We use an atomic measure that annihilates every polynomial,
an approach also used in
\cite[proof of Theorem~VI]{Koosis1964}.
We choose an entire function that vanishes at all but one of
its atoms, so that its integral against the measure is nonzero.
The main task is to control the function between the atoms while
keeping the weight increasing and its logarithmic integral finite.
Rescaled Chebyshev polynomials provide this control while
preserving zero exponential type. To obtain separation in the
weighted $L^p$ spaces, we smooth the measure inside intervals
on which the weight is constant. A bounded modification of the
weight then gives strict monotonicity.

\section{Proof of the main result}\label{sec:proof}

We use the following estimates for the Chebyshev factors. The bound on
complex disks will give both convergence and zero exponential type of
the infinite products.

\begin{lemma}[Chebyshev factors]\label{lem:block}
Let $A,a\in\R$ satisfy $2\leq A\leq a/2$, and let $m\geq1$
be an integer. Let $T_m$ denote the Chebyshev polynomial of the
first kind, characterized by
$T_m(\cos\theta)=\cos(m\theta)$ for $\theta\in\R$, and set
\begin{equation}\label{eq:block}
Q(z)=Q_{A,a,m}(z)=
\frac{T_m\!\left(\dfrac{4a^2+A^2-2z^2}{4a^2-A^2}\right)}
{T_m\!\left(\dfrac{4a^2+A^2}{4a^2-A^2}\right)},
\qquad z\in\C.
\end{equation}
Then $Q(0)=1$, and $Q$ has degree $2m$ and simple zeros
$\pm\xi_j$, $1\leq j\leq m$, with $A<\xi_j<2a$. Moreover, for real $x$,
\begin{equation}\label{eq:blockbounds}
\begin{aligned}
|Q(x)|&\leq1 &&(|x|\leq2a),\\
|Q(x)|&\leq2e^{-mA/a} &&(A\leq|x|\leq2a),
\end{aligned}
\end{equation}
and, for every $r\geq0$,
\begin{equation}\label{eq:complexblock}
\max_{|z|\leq r}|Q(z)|\leq e^{mr/a},
\qquad
\max_{|z|\leq r}|Q(z)-1|\leq e^{mr/a}-1.
\end{equation}
\end{lemma}

\begin{proof}
We also have $T_m(\cosh s)=\cosh(ms)$ for $s\geq0$. Write
$u(z)=(4a^2+A^2-2z^2)/(4a^2-A^2)$. Since $u(A)=1$ and
$u(2a)=-1$, the roots of $T_m$ give the stated zeros of $Q$.
Put
\[
\eta(r)=\operatorname{arcosh}
\frac{4a^2+A^2+2r^2}{4a^2-A^2}.
\]
Then
\[
\eta(0)=\log\frac{2a+A}{2a-A}\geq A/a,
\qquad T_m(u(0))=\cosh(m\eta(0)).
\]
On $|x|\leq A$, the argument $u(x)$ lies between $1$ and $u(0)$,
so $0<Q(x)\leq1$. On $A\leq|x|\leq2a$, the numerator has modulus
at most one. Since $\cosh(m\eta(0))\geq1$ and
$\cosh(m\eta(0))\geq e^{mA/a}/2$, this proves
\eqref{eq:blockbounds}.

Direct differentiation gives
\[
\eta'(r)=\frac{2r}{\sqrt{(r^2+A^2)(r^2+4a^2)}}\leq\frac1a.
\]
Thus $\eta(r)-\eta(0)\leq r/a$, and
$\cosh(s+t)\leq e^t\cosh s$ for $s,t\geq0$ yields
$Q(ir)\leq e^{mr/a}$. In the factorization
$Q(z)=\prod_{j=1}^m(1-z^2/\xi_j^2)$ the coefficients alternate
in sign and have nonnegative absolute values. Hence, for $|z|\leq r$,
$|Q(z)|\leq Q(ir)$ and $|Q(z)-1|\leq Q(ir)-1$,
which proves \eqref{eq:complexblock}.

\end{proof}

\begin{proof}[Proof of Theorem~\ref{thm:main}]
The proof has four steps.

\textbf{1. The entire function.}
Set
\begin{equation}\label{eq:parameters}
\begin{gathered}
a_1=1,\qquad a_2=e^{100},\qquad
a_{n+1}=\exp\left(\frac{a_n}{\log a_n}\right)\quad(n\geq2),\\
L_n=\log a_n,\qquad A_n=L_n^4,\qquad
m_n=\left\lceil\frac{a_n}{L_n^2}\right\rceil\quad(n\geq2),
\end{gathered}
\end{equation}
and let $Q_n=Q_{A_n,a_n,m_n}$ for $n\geq2$. Since
$L_{n+1}=e^{L_n}/L_n\geq2L_n$ for $L_n\geq100$, we have
$L_n\geq100\,2^{n-2}$. Also $e^{L_n}\geq2L_n^4$, so
\begin{equation}\label{eq:parambounds}
a_n\geq2A_n,\qquad
A_{n+1}=(a_n/L_n)^4\geq a_n^2,\qquad
a_{n+1}\geq4a_n.
\end{equation}
The last inequality holds also for $n=1$. The preceding lower bound
and the definition of $m_n$ give
\begin{equation}\label{eq:summability}
\sum_{n\geq2}\frac{n}{L_n}<\infty,\qquad
\frac1{L_n^2}\leq\frac{m_n}{a_n}
<\frac1{L_n^2}+\frac1{a_n},
\qquad \frac{m_nA_n}{a_n}\geq L_n^2.
\end{equation}

For $n\geq2$, write $B_n(z)=\prod_{j=1}^n(1-z^2/a_j^2)$ and
$H_{n-1}(z)=\prod_{j=2}^{n-1}Q_j(z)$, with empty products equal
to one. Their combined degree is
\[
D_n:=\deg(B_nH_{n-1})=2n+2\sum_{j=2}^{n-1}m_j.
\]
For $n\geq3$, the recurrence gives $a_{n-1}=L_nL_{n-1}$.
Consequently,
\[
m_{n-1}\sim\frac{L_n}{L_{n-1}},\qquad
\frac{m_n}{m_{n-1}}\sim
\frac{e^{L_n}L_{n-1}}{L_n^3}\longrightarrow\infty.
\]
The last term therefore dominates the sum defining $D_n$, and
$D_n\sim2m_{n-1}\sim2L_n/L_{n-1}$. Since
$\log L_n=L_{n-1}-\log L_{n-1}\sim L_{n-1}$, we obtain
\begin{equation}\label{eq:degreebalance}
\begin{aligned}
D_n\log(1+A_n)
&\sim\frac{2L_n}{L_{n-1}}\,4\log L_n\sim8L_n,\\
D_n\log(1+2a_n)
&\sim\frac{2L_n^2}{L_{n-1}}=o(L_n^2).
\end{aligned}
\end{equation}

These choices balance the requirements of the construction.
The identity $(n+1)L_{n+1}/a_n=(n+1)/L_n$ will make the
logarithmic cost of increasing the weight summable.
The summability of $m_n/a_n$ controls the exponential type
of the product of the Chebyshev factors. At the same time,
$m_nA_n/a_n\geq L_n^2$, whereas
$D_n\log(1+2a_n)=o(L_n^2)$, so the decay of $Q_n$ on
$[A_n,2a_n]$ dominates the growth of $B_nH_{n-1}$.

Define
\begin{equation}\label{eq:functions}
B(z)=\prod_{n\geq1}(1-z^2/a_n^2),\quad
H(z)=\prod_{n\geq2}Q_n(z),\quad
F(z)=\frac{B(z)H(z)}{z-1}.
\end{equation}
Equations \eqref{eq:complexblock} and \eqref{eq:summability} give
local uniform convergence of both products. The zeros of $B$ are
precisely the simple zeros $\pm a_n$, and the singularity at $1$
is removable:
\[
F(z)=-(z+1)\prod_{n\geq2}(1-z^2/a_n^2)Q_n(z).
\]
To check zero exponential type, fix an integer $N\geq2$ and split
off the factors with index at most $N$.
The remaining factors satisfy
\[
\left|\prod_{n>N}(1-z^2/a_n^2)Q_n(z)\right|
\leq\exp\left(
|z|\sum_{n>N}\left(\frac2{a_n}+\frac{m_n}{a_n}\right)\right),
\]
where we used $|1-z^2/a_n^2|\leq e^{2|z|/a_n}$.
The coefficient of $|z|$ tends to zero, while each initial product
is a polynomial. Thus $B,H,F$ have zero exponential type.
Moreover, $Q_n(1)>0$ and
$\sum_{n\geq2}|Q_n(1)-1|<\infty$, so $H(1)>0$.

\medskip
\noindent\textbf{2. The weight and the real-axis estimates.}
Put $w_n=a_n^n$ for $n\geq2$. Set $W=w_2$ on $[0,a_2+1]$.
For each $n\geq2$, put
\begin{equation}\label{eq:ramps}
W(x)=w_n+\frac{x-a_n-1}{a_n-1}(w_{n+1}-w_n),
\qquad a_n+1\leq x\leq2a_n,
\end{equation}
and set $W=w_{n+1}$ on $[2a_n,a_{n+1}+1]$. Extend $W$ evenly.
The heights increase, so $W$ is continuous, nondecreasing on
$[0,\infty)$, and at least one. It is constant within distance
$1/2$ of every $\pm a_n$, including $\pm1$. These neighborhoods
will allow us to smooth the separating measure.

For $2a_{n-1}\leq x\leq2a_n$, we have $W(x)\geq a_n^n$, whence
\[
\frac{\log W(x)}{\log x}\geq
\frac{nL_n}{L_n+\log2}\longrightarrow\infty.
\]
On $[a_n,a_{n+1}]$ we have $W\leq w_{n+1}$. Therefore
\begin{equation}\label{eq:logintegralbound}
\int_{a_2}^{\infty}\frac{\log W(x)}{1+x^2}\,dx
\leq\sum_{n\geq2}\frac{(n+1)L_{n+1}}{a_n}
=\sum_{n\geq2}\frac{n+1}{L_n}<\infty.
\end{equation}
Continuity on compact intervals and evenness prove
\eqref{eq:growth} and \eqref{eq:logint} on the whole line.

We next estimate $F$ on two regions. For $n\geq3$,
$2a_{n-1}\leq A_n\leq a_n/2$. For $2\leq|x|\leq2a_n$,
all factors with index greater than $n$ have modulus at most one.
The polynomial $B_nH_{n-1}$ is normalized at zero, has degree
$D_n$, and all its zeros are real and have modulus at least one.
Since $|x-1|\geq1$ for $|x|\geq2$, we obtain
\[
\left|\frac{B_n(x)H_{n-1}(x)}{x-1}\right|
\leq(1+|x|)^{D_n}.
\]
On the gap $2a_{n-1}\leq|x|\leq A_n$, use $|Q_n|\leq1$.
On $A_n\leq|x|\leq2a_n$, use the second bound in
\eqref{eq:blockbounds}. Equations \eqref{eq:summability} and
\eqref{eq:degreebalance} give, for all sufficiently large $n$,
\begin{equation}\label{eq:realbounds}
\begin{aligned}
|F(x)|&\leq(1+A_n)^{D_n}\leq a_n^9
&& (2a_{n-1}\leq|x|\leq A_n),\\
|F(x)|&\leq2e^{D_n\log(1+2a_n)-m_nA_n/a_n}
\leq2e^{-L_n^2/2}
&& (A_n\leq|x|\leq2a_n).
\end{aligned}
\end{equation}
On the first region $W=a_n^n$, so $|F|/W\leq a_n^{9-n}$;
on the second $W\geq1$. These bounds tend to zero and cover the
whole tail. The two regions, including both half-lines, have
length at most $2a_n$ and $4a_n$, respectively. Their contributions
to $\|F/W\|_1$ are therefore at most
\[
2a_n^{10-n}\qquad\text{and}\qquad8e^{L_n-L_n^2/2}.
\]
Both series converge: for $n\geq11$ the first is at most $2/a_n$,
and the second is at most $8e^{-L_n^2/4}$.
The remaining intervals form a compact set. We conclude that
$F/W\in C_0(\R)\cap L^1(\R)$.
Since $F/W$ is bounded, it also belongs to $L^p(\R)$ for every
$1\leq p<\infty$. Thus $F\in C_W^0$ and
$F\in L_W^p$ throughout this range.

\medskip
\noindent\textbf{3. A functional separating all the weighted spaces.}
Let $\Lambda=\{\pm a_n:n\geq1\}$, and let $\delta_t$ denote the
unit point mass at $t$. Differentiating the product for $B$ gives
\[
|B'(a_n)|=
\frac{2a_n^{2n-3}}{\prod_{j<n}a_j^2}
\prod_{j<n}\left(1-\frac{a_j^2}{a_n^2}\right)
\prod_{j>n}\left(1-\frac{a_n^2}{a_j^2}\right).
\]
By \eqref{eq:parambounds}, each of the two correction products is
at least $1-\sum_{k\geq1}16^{-k}=14/15$. Also
$2\sum_{j<n}\log a_j=O(L_{n-1})=o(L_n)$, so eventually
$\prod_{j<n}a_j^2\leq a_n$. Consequently,
\[
|B'(a_n)|\geq2(14/15)^2a_n^{2n-4}\geq a_n^{2n-4}.
\]
It follows that the real signed measure
\begin{equation}\label{eq:measure}
\mu=\sum_{t\in\Lambda}\frac{\delta_t}{B'(t)}
\quad\text{satisfies}\quad
\int_\R W\,d|\mu|<\infty,
\end{equation}
because its weighted variation has a tail bounded by
$2\sum_n a_n^{4-n}$. The growth of $W$ also implies that
every absolute moment of $\mu$ is finite.

For an integer $N\geq2$, Lagrange interpolation gives
\[
P(z)=\sum_{B_N(t)=0}P(t)
\frac{B_N(z)}{B_N'(t)(z-t)}
\qquad(\deg P\leq2N-1).
\]
Comparing coefficients of $z^{2N-1}$ when $\deg P\leq2N-2$ yields
\begin{equation}\label{eq:finiteidentity}
\sum_{B_N(t)=0}\frac{P(t)}{B_N'(t)}=0.
\end{equation}
For a fixed node $t=\pm a_k$ and $N\geq k$,
\[
\frac1{B_N'(t)}
=\frac{\prod_{j>N}(1-t^2/a_j^2)}{B'(t)}.
\]
The tail product belongs to $(0,1]$ and tends to one.
Extend the finite sums by zero at the remaining nodes.
They are dominated by the summable sequence
$|P(t)|/|B'(t)|\leq\|P\|_W W(t)/|B'(t)|$.
Dominated convergence in \eqref{eq:finiteidentity} gives
$\int P\,d\mu=0$ for every polynomial. On the other hand,
$F$ vanishes at every atom except $1$, where $F(1)=B'(1)H(1)$.
Thus
\begin{equation}\label{eq:pairing}
\int_\R F\,d\mu=H(1)>0.
\end{equation}

The measure $\mu$ already separates $F$ from the polynomials
in $C_W^0$. To obtain bounded separating functionals on the
spaces $L_W^p$, $1\leq p<\infty$, we smooth $\mu$.

Choose a nonnegative $\eta\in C_c^\infty(-1,1)$ with $\int\eta=1$, and set
\begin{equation}\label{eq:smoothing}
\eta_\varepsilon(x)=\varepsilon^{-1}\eta(x/\varepsilon),
\qquad g_\varepsilon=\mu*\eta_\varepsilon,
\qquad0<\varepsilon<1/2.
\end{equation}
Since $W$ is constant on each translated support,
\[
Wg_\varepsilon
=
\sum_{t\in\Lambda}\frac{W(t)}{B'(t)}
\eta_\varepsilon(\,\cdot-t).
\]
The triangle inequality therefore gives
\[
\|Wg_\varepsilon\|_q
\leq\|\eta_\varepsilon\|_q\int_\R W\,d|\mu|<\infty
\qquad(1\leq q\leq\infty).
\]
Finite absolute moments justify Fubini, and
\[
\int_\R P(x)g_\varepsilon(x)\,dx
=\int_\R\left(\int_\R
P(t+y)\eta_\varepsilon(y)\,dy\right)d\mu(t)=0,
\]
because the inner integral is a polynomial in $t$.
Continuity of $F$ and the summable bound
\[
\frac1{|B'(t)|}\int_\R|F(t+y)|\eta_\varepsilon(y)\,dy
\leq\|F\|_W\frac{W(t)}{|B'(t)|}
\]
give $\int Fg_\varepsilon\to H(1)$. Fix a sufficiently small
$\varepsilon$ and put $g=g_\varepsilon$. For $1/p+1/q=1$,
H\"older's inequality yields
\begin{equation}\label{eq:distance}
\inf_{P\in\C[x]}\left\|\frac{F-P}{W}\right\|_{L^p}
\geq\frac{\left|\int_\R Fg\,dx\right|}{\|Wg\|_q}>0
\qquad(1\leq p\leq\infty).
\end{equation}
The choice of $g$ is independent of $p$. For $p=\infty$, the
left-hand side is the weighted uniform distance, since the errors
are continuous.

\medskip
Set $\widetilde W(x)=(2-e^{-x^2})W(x)$. This weight is even,
continuous, and strictly increasing on $[0,\infty)$. Since
$W\leq\widetilde W\leq2W$, the corresponding weighted spaces
coincide with equivalent norms. Hence the membership and
positive-distance conclusions proved above are preserved.
The growth condition \eqref{eq:growth} is also preserved, and
\[
\int_\R\frac{\log\widetilde W(x)}{1+x^2}\,dx
\leq
\int_\R\frac{\log W(x)}{1+x^2}\,dx+\pi\log2
<\infty.
\]

\medskip
\noindent\textbf{4. Approximation by small frequencies.}
Let $X$ be either $C_{\widetilde W}^0$ or $L_{\widetilde W}^p$,
$1\leq p<\infty$. We use the weighted approximation characterization
stated in \cite[Introduction, p.~831]{Borichev2001}.
Its hypotheses are satisfied: $\widetilde W:\R\to[1,\infty)$
is even and continuous, satisfies \eqref{eq:growth} and
\eqref{eq:logint}, and $F$ is an entire function of zero
exponential type belonging to $X$. The characterization
therefore gives
\[
F\in
\bigcap_{a>0}
\overline{\operatorname{span}\{e^{itx}:t\in\R,\ |t|\leq a\}}^{\,X}.
\]
The approximation implication used here also follows by duality
from \cite[Theorem~I]{deBranges1959}.
Taking $\widetilde W$ as the weight in the theorem completes
the proof.
\end{proof}

\Needspace{12\baselineskip}
\begin{proof}[Proof of Corollary~\ref{cor:strictinclusion}]
For each integer $k\geq1$ and $h>0$,
$$
\left(\frac{e^{ihx}-1}{ih}\right)^k\longrightarrow x^k
\quad(h\to0),\qquad
\left|\frac{e^{ihx}-1}{ih}\right|^k\leq|x|^k.
$$
The convergence holds in $X$: use \eqref{eq:growth} to bound the
tails in $C_W^0$, and dominated convergence in $L_W^p$.
For a fixed $a>0$, take $h$ small enough that $kh\leq a$.
The approximants have frequencies $0,h,\ldots,kh$, and the constant
function is the exponential of frequency zero. Thus every polynomial
belongs to every exponential closure in the statement of the corollary. The function
$F$ belongs to their intersection and not to the polynomial closure.
\end{proof}

\section*{Acknowledgments}
P.I.~acknowledges partial support from NSF grant DMS-2554183, and a Humboldt Research Fellowship for Experienced Researchers. The authors acknowledge the use of AI tools in this project.

\Needspace{15\baselineskip}


\begin{thebibliography}{9}

\bibitem{Borichev2000}
A.~Borichev,
\emph{On weighted polynomial approximation with monotone weights},
Proc. Amer. Math. Soc. \textbf{128} (2000), no.~12, 3613--3619.

\bibitem{Borichev2001}
A.~Borichev,
\emph{On the closure of polynomials in weighted spaces of functions on the
real line},
Indiana Univ. Math. J. \textbf{50} (2001), no.~2, 829--846.

\bibitem{deBranges1959}
L.~de~Branges,
\emph{The $a$-local operator problem},
Canad. J. Math. \textbf{11} (1959), 583--592.

\bibitem{Koosis1964}
P.~Koosis,
\emph{Sur l'approximation pond\'er\'ee par des polynomes et par des
sommes d'exponentielles imaginaires},
Ann. Sci. \'Ecole Norm. Sup. (3) \textbf{81} (1964), no.~4,
387--408.

\end{thebibliography}
\end{document}